\documentclass[10pt,twoside,reqno]{amsart}
\usepackage{amssymb}
\calclayout

\numberwithin{equation}{section}
\makeatletter

\renewcommand{\@secnumfont}{\bfseries}

\renewcommand{\section}{\@startsection{section}{1}%
  {0mm}{.7\linespacing\@plus\linespacing}{.5\linespacing}
  {\normalfont\bfseries\centering}}

\newcommand{\bibsection}{\@startsection{section}{1}%
  {0mm}{.7\linespacing\@plus\linespacing}{.5\linespacing}
  {\normalfont\scshape\centering}}

\renewcommand{\@biblabel}[1]{#1.}

\newtheorem{thm}{\bf Theorem}[section]

\newtheorem{cor}[thm]{\bf Corollary}

\begin{document}

\vspace{1.3cm}

\title
        {$p$-adic invariant integral on $ \mathbb{Z}_p$ associated with the Changhee's $q$-Bernoulli polynomials}

\author{J.J. Seo, T.Kim}

\thanks{\scriptsize }

\address{1\\Department of Applied Mathematics\\
            Pukyong National University\\
            Busan 608-737, Republic of Korea}
\email{seo2011@pknu.ac.kr}

\address{2\\Department of  Mathematics\\
            Kwangwoon National University\\
           Seoul 139-701, Republic of Korea}
\email{tkkim@kw.ac.kr}


\maketitle

\begin{abstract} In this paper, we study some properties of Changhee's $q$-Bernou\\lli polynomials which are derived from $p$-adic invariant integral on $\mathbb{Z}_p$. By using these properties, we give some interesting identities related to higher-order $q$-Bernoulli polynomials.
\end{abstract}

\pagestyle{myheadings}
\markboth{\centerline{\scriptsize J.J. Seo, T.Kim}}
          {\centerline{\scriptsize $p$-adic invariant integral on $ \mathbb{Z}_p$ associated with the Changhee's $q$-Bernoulli polynomials}}
\bigskip
\bigskip
\medskip
\section{\bf Introduction}
\bigskip
\medskip

Throughout this paper, $ \mathbb{Z}_p, \mathbb{Q}_p, \mathbb{C}$ and $\mathbb{C}_p$ will respectively denote the ring of $p$-adic integers, the field of $p$-adic numbers, the complex number field and the completion of algebraic closure of $\mathbb{Q}_p$. Let $\nu_p$ be the normalized exponential valuation of $\mathbb{C}_p$ with $|p|_p=p^{^{-\nu_p (p)}}=1/p.$ When one talks of $q$-extension, $q$ is variously considered as an indeterminate, a complex number $q\in \mathbb{C}$, or a $p$-adic number $q\in\mathbb{C}_p$. If $q\in\mathbb{C}$, one normally assumes $|q|<1$. If  $q\in \mathbb{C}_p$, one normally assumes $|1-q|_p<p^{-1/p-1},$ so that $g^x=exp(xlogq)$ for $x \in\mathbb{Z}_p $.
We use the following notation in this paper:\
\begin{equation*}[x]_q=[x:q]= {{1-q^x}\over{1-q}}, \ \ (\textnormal{see}\ [8]-[16]).
\end{equation*}
Note that $\lim_{q\rightarrow 1} {[x]_q=x}$ for any $x$ with $|x|_p \leq 1$ in the present $p$-adic case.\
Let $d$ be a fixed integer and let $p$ be a fixed prime number. We set
\begin{equation*}
\mathbb{X}_d=\lim_{\overleftarrow{N}} {\mathbb{Z}\diagup{dp^{N}\mathbb{Z},}}
\ \ \  \mathbb{X^{*}}={\bigcup_{\substack{0< a< dp\\  (a,p)=1}}} {\ (a+dp\mathbb{Z}_p)},\\
\end{equation*}
\begin{equation*}
a+dp^{N}\mathbb{Z}_p=\left\{x\in X | x\equiv a (mod\ dp^{N}) \right\},\\
\end{equation*}
where $a\in \mathbb{Z}$ lies in  $0\leq a < {dp^{N}}$,\ \ (\textnormal{see}\ [1]-[11]).
Let $UD(\mathbb{Z}_p)$ be the space of uniformly differentiable functions on $\mathbb{Z}_p$.\
For $f\in UD(\mathbb{Z}_p)$, the $p$-adic invariant integral on  $\mathbb{Z}_p$ is defined by\
\begin{equation}
I_0(f)=\int_{\mathbb{Z}_p} f(x) d\mu_0 (x) = \lim_{N\rightarrow \infty} {1\over {p^N}} \sum_{x=0}^{p^N-1} f(x)= \int_{\mathbb{X}} f(x) d\mu_0 (x).
\end{equation}
\ Thus, by$(1.1)$, we see \\
\begin{equation}
I_0(f_1)= I_0(f)+f^{\prime}(0),
\end{equation}
where $f_1(x)=f(x+1)$.
From $(1.2)$, we can derive\\
\begin{equation}
I_0(f_n)=  I_0(f)+\sum_{i=0}^{n-1}f^{\prime}(i),(\textnormal{see}\ [9],[10],[14],[15])
\end{equation}
where $f_n(x)=f(x+n),$ \ $f^{\prime}(i)={{df(x)}\over{dx}}|_{x=i}$.\\
As is well known, the Bernoulli polynomials of order $r$ are defined by the generating function to be
\begin{equation}
\left({{t}\over{e^{t}-1}} \right)^{r} e^{xt}=\sum_{n=0}^{\infty}B_n^{(r)}(x) {{t^{n}}\over{n!}}, \ (r\in\mathbb{Z}_{\geq 0} )
\end{equation}
In the special case, $x=0$,\ $B_n^{(r)}(0)=B_n^{(r)}$ are called the $n-$th Bernoulli numbers of order $r$(\textnormal{see}\ [1]-[16]).
By $(1.2)$, we easily set\\
\begin{equation}\begin{split}
&\underbrace{\int_{\mathbb{Z}_p}\cdots \int_{\mathbb{Z}_p}}_{r-times} e^{(x_1+ \cdots+x_r)t} \ d\mu_0 (x_1) \cdots d\mu_0 (x_r)\\
&\ \ \ \ \ \ \ \ \ \ \ \ \ \ \ \ \ \ \ \ \ \ =\left({{t}\over{e^{t}-1}} \right)^{r}=\sum_{n=0}^{\infty}B_n^{(r)} {{t^{n}}\over{n!}}.
\end{split}\end{equation}
From $(1.5)$, we have\\
\begin{equation}\begin{split}
&\underbrace{\int_{\mathbb{Z}_p}\cdots \int_{\mathbb{Z}_p}}_{r-times} {(x_1+ \cdots+x_r)^{n}} \ d\mu_0 (x_1) \cdots d\mu_0 (x_r)=B_n^{(r)}, (n\geq 0),
\end{split}\end{equation}
and\\
\begin{equation}\begin{split}
&\underbrace{\int_{\mathbb{Z}_p}\cdots \int_{\mathbb{Z}_p}}_{r-times} {(x_1+ \cdots+x_r+x)^{n}} \ d\mu_0 (x_1) \cdots d\mu_0 (x_r)=B_n^{(r)}(x), (n\geq 0).
\end{split}\end{equation}
In this paper, we consider Changhee $q$-Bernoulli number and polynomials which are derived from multivariate $p-$adic invariant integral on $\mathbb{Z}_p$.
 and investigate some properties of Changhee $q$-Bernoulli polynomials.

\section{\bf Changhee $q$-Bernoulli polynomials}
\bigskip
\medskip
Let $a_1 ,a_2 ,\cdots, a_k, b_1 ,b_2 , \cdots, b_k$ be positive integers. For $w\in \mathbb{Z}_p,$ we consider the modified Changhee $q$-Bernoulli polynomials which are derived from multivariate $p$-adic invariant integral on $\mathbb{Z}_p$ as follows:
\begin{equation}\begin{split}
&B_{n,q}^{(k)} \left(w  |  a_1, a_2 , \cdots, a_k ; b_1 ,b_2 , \cdots, b_k \right)\\
&=\underbrace{\int_{\mathbb{Z}_p}\cdots \int_{\mathbb{Z}_p}}_{k-times} q^{b_1 x_1 +\cdots+b_k x_k} \left[w+a_1 x_1 +\cdots +a_k x_k \right]_q^n d\mu_0 (x_1) \cdots d\mu_0 (x_k),
\end{split}\end{equation}
and
\begin{equation}\begin{split}
&B_{n,q}^{(k)} \left(  a_1, a_2 , \cdots, a_k ; b_1 ,b_2 , \cdots, b_k \right)\\
&=\underbrace{\int_{\mathbb{Z}_p}\cdots \int_{\mathbb{Z}_p}}_{k-times} q^{b_1 x_1 +\cdots+b_k x_k} \left[a_1 x_1 +\cdots +a_k x_k \right]_q^n d\mu_0 (x_1) \cdots d\mu_0 (x_k).
\end{split}\end{equation}
From (1.1) and(2.1), we have
\begin{equation}\begin{split}
&B_{n,q}^{(k)} \left(w  |  a_1, a_2 , \cdots, a_k ; b_1 ,b_2 , \cdots, b_k \right)\\
&=\underbrace{\int_{\mathbb{Z}_p}\cdots \int_{\mathbb{Z}_p}}_{k-times} q^{b_1 x_1 +\cdots+b_k x_k} \left[w+a_1 x_1 +\cdots +a_k x_k \right]_q^n d\mu_0 (x_1) \cdots d\mu_0 (x_k)\\
&= {1 \over(1-q)^n } \sum_{r=0}^n \binom{n}{r} (-q^w)^r \lim_{N\rightarrow \infty} \left(1\over P^N \right)^k\\
 & \ \ \ \ \ \ \ \ \ \ \ \ \ \ \ \ \ \ \ \ \ \times \sum_{x_1 ,x_2 ,\cdots, x_k=0}^{P^N-1} q^{ra_1x_1 +\cdots+ra_k x_k +b_1x_1 +\cdots+b_kx_k}\\
&= {1 \over(1-q)^n } \sum_{r=0}^n \binom{n}{r} (-q^w)^r \left({{\log q}\over{q-1}} \right)^k \prod_{j=1}^{k} \left({{ra_j+b_j}\over{[ra_j+b_j]_q}} \right)
\end{split}\end{equation}
Therefore, by (2.3), we obtain the following theorem.
\begin{thm}\label{THEOREM 2.1} For $a_1 ,a_2 , \cdots, a_k , b_1 ,b_2 , \cdots, b_k \in \mathbb{Z}_{\geq 0}$ and  $ w\in \mathbb{Z}_p$, we have
\begin{equation*}\begin{split}
&B_{n,q}^{(k)} (w \ | \  a_1 ,a_2 ,\cdots, a_k ; b_1 ,b_2 ,\cdots, b_k )\\
&=\left({\log q}\over {q-1} \right)^k {1\over {(1-q)^n}}\sum_{r=0}^{n}\binom{n}{r} (-q^w)^r \prod_{j=1}^{k} \left({{ra_j+b_j}\over{[ra_j+b_j]_q}} \right),
\end{split}\end{equation*}
where $n\geq0 , k\geq 1$.
\end{thm}
\bigskip
\medskip
\bigskip

Let us take $b_1=1, b_2 =2, \cdots, b_k =k$. Then, by theorem 2.1, we get
\begin{equation}\begin{split}
&B_{n,q}^{(k)} (w \ | \underbrace{\ 1, 1, \cdots, 1}_{k-times} ;1, 2,\cdots, k )\\
&=\left({\log q}\over {q-1} \right)^k {1\over {(1-q)^n}}\sum_{r=0}^{n}\binom{n}{r} (-q^w)^r \prod_{j=1}^{k} \left({{r+j}\over{[r+j]_q}} \right)\\
&=\left({\log q}\over {q-1} \right)^k {1\over {(1-q)^n}}\sum_{r=0}^{n}\binom{n}{r} (-q^w)^r {{\binom{r+k}{k}}\over{\binom{r+k}{k}_q}} {{k!}\over{[k]_q!}},
\end{split}\end{equation}
where ${\binom{r}{k}}_q={{[r]_q \cdots [r-k+1]_q}\over{[k]_q!}},$ $[k]_q!=[k]_q [k-1]_q \cdots [2]_q [1]_q$.
Therefore, by (2.4), we obtain the following corollary.
\begin{cor}\label{COROLLARY 2.2} For $n\geq 0, k\geq1 $, we have
\begin{equation*}
B_{n,q}^{(k)} (w \ |\underbrace{ \  1, 1, \cdots, 1}_{k-times} ; 1, 2, \cdots, k)
=\left({\log q}\over {q-1} \right)^k {1\over {(1-q)^n}}\sum_{r=0}^{n}\binom{n}{r} (-q^w)^r {{\binom{r+k}{k}}\over{\binom{r+k}{k}}_q} {{k!}\over{[k]_q!}}.
\end{equation*}
\end{cor}
The multiple Barnes' Bernoulli polynomials are defined by the generating function to be
\begin{equation}
\prod_{j=1}^{k} \left({w_j \over {e^{w_j t}-1}} \right) t^k e^{xt} = \sum_{n=0}^{\infty} B_n^{(k)} (x | w_1 ,\cdots, w_k) {t^n \over n!},
\end{equation}
where $w_i>0,$ $0<t<1$ (see [8],[9],[15])).
The numbers $B_n^{(k)} (0| w_1 ,w_2 ,\cdots, w_k)=B_n^{(k)}(w_1 ,\cdots, w_k)$ are called the Barnes' Bernoulli numbers of order $k$.
If we take $x=rx, w_k =r+k$ and $t=\log q$ in (2.5), then we have
\begin{equation*}
{{(r+1)\cdots (r+k)}\over{(q^{r+1}-1)\cdots (q^{r+k}-1)}} q^{rx} = \sum_{s=0}^{\infty} B_s^{(k)} (rx | r+1, \cdots, r+k ) {{(\log q)^s}\over{s!}}.
\end{equation*}
Note that
\begin{equation*}
\lim_{q\rightarrow 1} B_{n,q}^{(k)} (w | a_1 ,a_2 , \cdots, a_k ; b_1 ,b_2 , \cdots, b_k)=B_n^{(k)}(a_1, \cdots, a_k)
\end{equation*}
For $k=1$ and $w=0$, let $B_{n,q}=B_{n,q}^{(1)}(0 | 1).$ Then we have
\begin{equation*}
B_{n,q}=B_{n,q}^{(1)} (0 | 1) =\left({\log q \over 1-q} \right) n \sum_{m=0}^{\infty} q^m [m]_q^{n-1} -\log q\sum_{m=0}^{\infty}q^m [m]_q^n.
\end{equation*}
From $(2.1)$, we can derive
\begin{equation}\begin{split}
&\underbrace{\int_{\mathbb{Z}_p} \cdots \int_{\mathbb{Z}_p}}_{k-times} q^{b_1 x_1 + \cdots +b_k x_k} [a_1 x_1 +\cdots + a_k x_k]_q^n d\mu_0 (x_1) \cdots d\mu_0 (x_k)\\
&=(q-1)\underbrace{\int_{\mathbb{Z}_p} \cdots \int_{\mathbb{Z}_p}}_{k-times} q^{(b_1-a_1) x_1 + \cdots +(b_k-a_k) x_k}\\
 & \ \ \ \ \ \ \ \ \ \ \ \ \ \ \ \ \ \ \ \ \times[a_1 x_1 +\cdots + a_k x_k]_q^{n+1} d\mu_0 (x_1) \cdots d\mu_0 (x_k)\\\\
&\ \ \ \ \ + \underbrace{\int_{\mathbb{Z}_p} \cdots \int_{\mathbb{Z}_p}}_{k-times} q^{\sum_{j=1}^{k}(b_j-a_j)x_j} [a_1 x_1 +\cdots + a_k x_k]_q^n d\mu_0 (x_1) \cdots d\mu_0 (x_k)\\
&=(q-1)B_{n+1,q}^{(k)} (a_1 , \cdots, a_k ; b_1-a_1, \cdots, b_k-a_k)\\
&\ \ \ \ \ +B_{n,q}^{(k)} (a_1, \cdots, a_k ; b_1 -a_1, \cdots, b_k -a_k ).
\end{split}\end{equation}
Therefore, by (2.6), we obtain the following theorem.
\begin{thm}\label{THEOREM 2.3} For $n\geq 0, k\geq 1$, we have
\begin{equation*}\begin{split}
&B_{n,q}^{(k)} (a_1 ,\cdots, a_k ; b_1 ,\cdots, b_k )=(q-1)B_{n+1,q}^{(k)}(a_1, \cdots, a_k ; b_1-a_1 ,\cdots, b_k-a_k)\\
& \ \ \ \ \ \ \ \ \  \ \ \ \ \ \ \ \ \ \ \ \ \ \ \ \ \ \ \ \ \ \ \ \ \ \ \ \ \ +B_{n,q}^{(k)}(a_1,\cdots, a_k;b_1-a_1,\cdots, b_k-a_k).
\end{split}\end{equation*}
\end{thm}
It is easy to show that
\begin{equation}\begin{split}
&\sum_{j=0}^{i}\binom{i}{j}(q-1)^j \underbrace{\int_{\mathbb{Z}_p}\cdots\int_{\mathbb{Z}_p}}_{k-times} [a_1 x_1+\cdots a_k x_k]_q^{n-i+j} \\
  & \ \ \ \ \ \ \ \ \ \ \ \ \ \ \ \ \ \ \ \ \ \ \ \ \ \ \ \ \ \ \ \ \ \ \ \ \ \ \times q^{b_1 x_1+\cdot +b_k x_k} d \mu_0(x_1)\cdots d\mu_0(x_k)\\
&=\sum_{j=0}^{i-1}\binom{i-1}{j}(q-1)^j \underbrace{\int_{\mathbb{Z}_p}\cdots\int_{\mathbb{Z}_p}}_{k-times}[a_1 x_1+\cdots a_k x_k]_q^{n-i+j} \\
& \ \ \ \ \ \ \ \ \ \ \ \ \ \ \ \ \ \ \ \ \ \ \ \ \ \ \ \ \ \ \ \ \ \ \ \ \ \ \ \times q^{(b_1+a_1)x_1 +\cdots +(b_k+a_k)x_k} d \mu_0(x_1)\cdots d\mu_0(x_k)\\
&=\sum_{j=0}^{i-1}\binom{i-1}{j}(q-1)^jB_{n-i+j,q}^{(k)}(a_1,\cdots, a_k ; b_1+a_1 , \cdots, b_k+a_k).
\end{split}\end{equation}
Therefore, by (2.7), we obtain the following theorem
\begin{thm}\label{THEOREM2.4} For $i \geq 1, n \geq 0, k \geq 1$, we have
\begin{equation*}\begin{split}
&\sum_{j=0}^{i} \binom{i}{j}(q-1)^jB_{n-i+j, q}^{(k)}(a_1,\cdots, a_k ; b_1 , \cdots, b_k)\\
&=\sum_{j=0}^{i-1}\binom{i-1}{j}(q-1)^j B_{n-i+j, q} ^{(k)}(a_1 ,\cdots, a_k ; b_1+a_1 , \cdots, b_k+a_k).
\end{split}\end{equation*}
In the special case, $k=1$, we have
\begin{equation*}\begin{split}
\sum_{j=0}^n \binom{n}{j}(q-1)^j B_{j, q}^{(1)} (a_1,b_1)&=\sum_{j=0}^n \binom{n}{j}(q-1)^j \int_{\mathbb{Z}_p}[a_1 x]_q^j q^{b_1 x} d\mu_0 (x)\\
&=\int_{\mathbb{Z}_p} q^{(na_1+b_1)x} d\mu_0(x)\\
&=\left({{q-1}\over{\log q}} \right) {{na_1+b_1}\over{[na_1+b_1]}}.
\end{split}\end{equation*}
\end{thm}
From (2.1) and (2.3), we have
\begin{equation}\begin{split}
&B_{n,q}^{(k)} (w | a_1,\cdots, a_k ; b_1 , \cdots, b_k)\\
&=\left({{\log q}\over{q-1}} \right)^k {1\over (1-q)^n} \sum_{r=0}^n \binom{n}{r}(-q^w)^r \prod_{j=1}^k \left({{ra_j+b_j}\over{[ra_j+b_j]_q}} \right)\\
&=\sum_{i=0}^{n} \binom{n}{i}[w]_q^{n-i}q^{wi}  B_{i, q}^{(k)}(a_1,\cdots, a_k ; b_1 , \cdots, b_k)
\end{split}\end{equation}
and
\begin{equation}\begin{split}
&B_{n, q}^{(k)}(w | a_1 ,\cdots, a_k ; b_1 , \cdots, b_k)\\
&= \underbrace{\int_{\mathbb{X}_d}\cdots \int_{\mathbb{X}_d}}_{k-times}  q^{b_1 x_1 + \cdots + b_k x_k} [a_1 x_1 + \cdots + a_k x_k]_q^n d\mu_0 (x_1)\cdots d\mu_0(x_k)\\
&=[l]_q^{n-k} \sum_{i_1 , \cdots, i_k =0}^{l-1} q^{b_1 i_1 + \cdots +b_k i_k}\\
& \ \ \ \ \ \ \times \underbrace{\int_{\mathbb{Z}_p}\cdots \int_{\mathbb{Z}_p}}_{k-times}\left[{{{w+a_1 i_1+\cdots+a_k i_k}\over{l}}+a_1x_1+\cdots+a_kx_k} \right]_{ql}^n \\
  & \ \ \ \ \ \ \ \ \ \ \ \ \ \ \ \ \ \ \ \ \ \ \ \ \ \ \ \ \ \ \ \ \ \ \ \ \ \ \ \times q^{lb_1 x_1 + \cdots + lb_k x_k} d\mu_0 (x_1)\cdots d\mu_0(x_k)\\
&=[l]_q^{n-k} \sum_{i_1 , \cdots, i_k=0}^{l-1} q^{b_1 i_1 + \cdots +b_k i_k}\\
& \ \ \ \ \ \times B_{n, q^l}^{(k)} \left( {{w+a_1 i_1+\cdots+a_ki_k}\over l}  \  | \  a_1, \cdots, a_k; b_1 , \cdots, b_k \right)
\end{split}\end{equation}
(Distribution relation for Changhee $q$-Bernoulli polynomials)
An obvious generating function $F_q (w,t)$ of $B_{n,q}^{(k)}(w |a_1 , \cdots, a_k ; b_1 , \cdots, b_k)$ is obtained from Theorem 2.1 as follows;
\begin{equation}\begin{split}
F_q(w,t)&=\left( {{\log q}\over {q-1}}\right)^k \sum_{n=0}^\infty (q-1)^n B_{n,q}^{(k)} (w | a_1 , \cdots, a_k ; b_1 , \cdots, b_k) {{t^n}\over{n!}}\\
&=\left({{\log q}\over {q-1}} \right)^k e^{-t} \sum_{i=0}^\infty \left(\prod_{j=1}^k {{ia_j+b_j}\over {[ia_j+b_j]_q}} \right) q^{wi} {{t^i}\over {i!}}.
\end{split}\end{equation}
Note that
\begin{equation*}
\lim_{q \rightarrow 1}F_q (w, t) ={{a_1 a_2 \cdots a_k t^k}\over{(e^{a_1 t}-1) \cdots (e^{a_k t} -1)}}.
\end{equation*}
Differentiating both sides with respect to $t$ in $(2.10)$ and comparing coefficients, we have
\begin{equation}\begin{split}
&q^w B_{n, q}^{(k)} (w | a_1 , \cdots, a_k ; b_1 , \cdots, b_k )- B_{n, q}^{(k)} (w | a_1 ,a_2 , \cdots, a_k ; b_1 -a_1 , \cdots, b_k -a_k )\\
&=(q-1) B_{n+1, q}^{(k)} (w | a_1 , \cdots, a_k ; b_1 -a_1 , \cdots, b_k -a_k ).
\end{split}\end{equation}
It is easy to show that
\begin{equation}
\sum_{i=0}^{n}\binom{n}{i} (q-1)^i B_{i, q}^{(k)} (a_1, \cdots, a_k ; b_1, \cdots, b_k ) = \left({{\log q}\over{q-1}} \right)^k \prod_{j=1}^k \left({{na_j+b_j}\over{[na_j +b_j]_q}} \right).
\end{equation}
Therefore, by (2.11) and (2.12), we obtain the following theorem.
\begin{thm}\label{THEOREM 2.5} For $n \geq 0, k \geq 1$, we have
\begin{equation*}\begin{split}
&q^w B_{n, q}^{(k)}(w | a_1 , \cdots, a_k ; b_1 , \cdots, b_k )- B_{n, q} ^{(k)} (w | a_1 , \cdots, a_k ; b_1 -a_1 , \cdots, b_k -a_k )\\
&= (q-1) B_{n+1, q}^{(k)} (w | a_1 , \cdots, a_k ; b_1 -a_1 , \cdots, b_k -a_k)\\
\end{split}\end{equation*}
Moreover,
\begin{equation*}
\sum_{i=0}^n \binom{n}{i}(q-1)^i B_{i, q}^{(k)} (a_1 , \cdots, a_k ; b_1 , \cdots, b_k )= \left({{\log q}\over{q-1}} \right)^k \prod_{j=1}^k \left({{na_j+b_j}\over{[na_j+b_j]_q}} \right)
\end{equation*}
\end{thm}
Let us assume that $a_1=\cdots=a_k=1, b_1 =h, b_2 =h-1, \cdots, b_k =h-k+1$, where $h \in \mathbb{Z}$. Then we have
\begin{equation}\begin{split}
&B_{n, q}^{(k)}(w | \underbrace{1, \cdots, 1}_{k-times} ; h, h-1, \cdots, h-k+1)\\
&=\left({{\log q}\over{q-1}}\right)^k {1\over{(1-q)^n}} \sum_{r=0}^n \binom{n}{r} (-q^w)^r \prod_{j=1}^k \left( {{r+h-j+1}\over{[r+h-j+1]_q}}\right)\\
&=\left({{\log q}\over{q-1}}\right)^k {{1}\over{(1-q)^n}}\sum_{r=0}^n \binom{n}{r} (-q^w)^r {{\binom{r+h}{k}}\over{\binom{r+h}{k}_q}} {{k!}\over{[k]_q!}}.
\end{split}\end{equation}
Therefore, by (2.13), we obtain the following theorem.
\begin{thm}\label{THEOREM 2.6} For $n \in \mathbb{Z} , k \geq 1, n\geq 0$, we have
\begin{equation*}\begin{split}
&B_{n, q}^{(k)}(w | 1, 1, \cdots, 1 ; h, h-1, \cdots, h-k+1 )\\
&= \left({{\log q}\over{q-1}}\right)^k {1 \over {(1-q)^n}} \sum_{r=0}^n \binom{n}{r}(-q^w)^r {{\binom{r+h}{k}}\over{\binom{r+h}{k}_q}}{{k!}\over{[k]_q!}}.
\end{split}\end{equation*}
\end{thm}
Let us use the notation of $q$-product as follows :
\begin{equation}
(a;q)_n = (1-a)(1-aq)\cdots (1-aq^{n-1})=\prod_{i=0}^{n-1}(1-aq^i).
\end{equation}
From (2.13), we have
\begin{equation}\begin{split}
&B_{n, q}^{(k)}(w | \underbrace{1, 1, \cdots, 1}_{k-times} ; h, h-1, \cdots, h-k+1)\\
&=(\log q)^k {1\over{(1-q)^n}}\sum_{r=0}^n \binom{n}{r} q^{wr} (-1)^{r+k} {(r+h)_k\over{(q^{r+h-k+1}:q)_k}}   \\
&=(\log q)^k {1\over{(1-q)^n}}\sum_{r=0}^n \binom{n}{r} q^{wr} (-1)^{r+k} k! \binom{r+h}{k} \\
&\ \ \ \ \ \ \ \ \ \ \ \ \ \ \ \ \ \ \ \ \ \ \ \ \ \ \ \ \ \ \ \ \ \ \ \ \ \ \ \ \ \times \sum_{m=0}^\infty \binom{m+k-1}{m}_q q^{m(r+h-k+1)}\\
&=k!(\log q)^k {1\over{(1-q)^n}}\sum_{m=0}^\infty \binom{m+k-1}{m}_q q^{m(h-k+1)} \\
&\ \ \ \ \ \ \ \ \ \ \ \ \ \ \ \ \ \ \ \ \ \ \ \ \ \ \ \ \ \ \ \ \ \ \ \ \ \ \ \  \times \sum_{r=0}^{n} \binom{n}{r} (-1)^{r+k} \binom{r+h}{k} q^{(m+w)r}.
\end{split}\end{equation}
By (2.15), we get
\begin{equation}\begin{split}
&B_{n, q}^{(k)} (w | \underbrace{1, 1, \cdots, 1}_{k-times} ; h, h-1, \cdots, h-k+1)\\
&={{(\log q)^k}\over {(1-q)^n}}k! \sum_{m=0}^{\infty} \binom{m+k-1}{m} q^{m(h-k+1)} \sum_{r=0}^n \binom{n}{r}\binom{r+h}{k} (-1)^{r+k} q^{(m+w)r}.
\end{split}\end{equation}
Therefore, by (2.16), we obtain the following theorem.
\begin{thm}\label{THEOREM 2.7} For $n \geq 0 , k \geq 1$, we have
\begin{equation}\begin{split}
&B_{n, q}^{(k)}(w | 1, 1, \cdots, 1 ; h, h-1, \cdots, h-k+1 )\\
&={{(\log q)^k}\over{(1-q)^m}} k! \sum_{m=0}^\infty \binom{m+k-1}{m}_q q^{m(h-k+1)} \sum_{r=0}^n \binom{n}{r} \binom{r+h}{k}(-1)^{r+h} q^{(m+w)r}.
\end{split}\end{equation}
\end{thm}
\hskip -1pc {\bf Remark.} From (2.15), we can derive
\begin{equation}\begin{split}
&q^r B_{n, q}^{(k)} (w+1 | \underbrace{1, 1, \cdots, 1}_{k-times} ; h, h-1, \cdots, h-k+1)\\
&={{(\log q)^k}\over{(1-q)^n}}q^h \sum_{j=0}^n \binom{n}{j} q^{(w+1)j} (-1)^{j+k} {{(j+h)_k}\over{(q^{j+h} ; q^{-1})_k}}\\
&={{(\log q)^k}\over{(1-q)^n}}\sum_{j=0}^n \binom{n}{j} (-1)^{j+k} q^{wj} \left([j+h]_q (q-1)+1 \right) {{(j+h)_k}\over{(q^{j+h} ; q^{-1})_k}}\\
&={{(q-1)(\log q)^k}\over{(1-q)^n}}\sum_{j=0}^n \binom{n}{j} (-1)^{j+k}q^{wj} {{(j+h-1)_{k-1}}\over{(q^{j+h-1} ; q^{-1} )_{k-1}}} (j+h)\\
&+{{(\log q)^k}\over{(1-q)^n}} \sum_{j=0}^n \binom{n}{j}(-1)^{j+k}q^{wj}{{(j+h)_k}\over{(q^{j+h} ; q^{-1})_k}} \\
&=h(1-q)\log q B_{n ,q}^{(k-1)} (w | \underbrace{1, 1, \cdots, 1}_{k-times} ; h-1, h-2, \cdots, h-k)\\
&-n \log q B_{n-1, q}^{(k-1)} (w | \underbrace{1, 1, \cdots, 1}_{k-times} ; h, h-1, \cdots, h-k+1)\\
&+B_{n, q}^{(k)}(w | \underbrace{1, 1, \cdots, 1}_{k-times} ; h, h-1, \cdots, h-k+1).
\end{split}\end{equation}

\end{document}